\documentclass[preprint,11pt,authoryear,3p]{elsarticle}
\usepackage[font=small,labelfont=bf]{caption}
\usepackage[font=small]{subcaption}

\usepackage{amsthm} 
\usepackage{algorithm, algpseudocode, amsmath}
\usepackage[english]{babel}
\usepackage[utf8x]{inputenc}

\usepackage[T1]{fontenc}
\usepackage{graphicx}

\usepackage{xcolor} 
\usepackage[english]{babel}
\usepackage{nomencl} 
\usepackage{amsmath} 
\usepackage{booktabs} 
\usepackage[]{comment} 
\usepackage{amssymb}
\usepackage{verbatim}

\usepackage{lscape}
\usepackage{placeins}

\usepackage{mathrsfs}

\usepackage{lipsum}

\usepackage{multicol}
\usepackage{pgfplots}
\usepgfplotslibrary{fillbetween}
\pgfplotsset{
   compat=newest,
   xlabel near ticks,
   ylabel near ticks
}

\makenomenclature

\usepackage{setspace}
\usepackage{booktabs}
\usepackage{multirow}
\usepackage{siunitx}

\usepackage{adjustbox}

\usepackage{graphicx}
\usepackage[colorinlistoftodos]{todonotes}
\usepackage[colorlinks=true, allcolors=blue]{hyperref}
\usepackage{tikz}  
\usetikzlibrary{calc} 
\usepackage{circuitikz} 
\usepackage{breqn} 
\usepackage{array}
\usepackage{bm}
\usepackage{nomencl} 
\makenomenclature

\usepackage{etoolbox}
\renewcommand\nomgroup[1]{%
  \item[\bfseries
  \ifstrequal{#1}{A}{\textit{Sets}}{%
  \ifstrequal{#1}{B}{\textit{Parameters}}{%
  \ifstrequal{#1}{E}{\textit{Variables}}{
  \ifstrequal{#1}{C}{\textit{Uncertain Parameters}}{%
  \ifstrequal{#1}{D}{\textit{none}}{%
  {}}}}}}%
]}

\newtheorem{prop}{Proposition}
\newtheorem{definition}{Definition}
\newtheorem{remark}{Remark}

\usepackage{todonotes}

\usepackage{interval}

\journal{European Journal of Operational Research}
\makeatletter
\newcommand\deletedfootnote[1]{%
    \ifthenelse{\boolean{Changes@optiondraft}}
    {\footnote{\deleted{#1}}}
    {}}
\makeatother

\begin{document}
\begin{frontmatter}

\title {Exponential Conic Relaxations for Signomial Geometric Programming} 

\author[label1]{Milad Dehghani Filabadi}
\address[label1]{Department of Integrated Systems Engineering, The Ohio State University, Columbus, OH, USA}
\ead{dehghanifilabadi.1@osu.edu}

\author[label1]{Chen Chen\corref{cor1}}
\cortext[cor1]{Corresponding author}
\ead{chen.8018@osu.edu}

\begin{abstract}
Signomial geometric programming (SGP) is a computationally challenging, NP-Hard class of nonconvex nonlinear optimization problems. SGP can be solved iteratively using a sequence of convex relaxations; consequently, the strength of such relaxations is an important factor to this iterative approach. Motivated by recent advances in solving exponential conic programming (ECP) problems, this paper develops a novel convex relaxation for SGP. Unlike existing work on relaxations, the base model in this paper does not assume bounded variables. However, bounded variables or monomial terms can be used to strengthen the relaxation by means of additional valid linear inequalities. We show how to embed the ECP relaxation in an iterative algorithm for SGP; leveraging recent advances in interior point method solvers, our computational experiments demonstrate the practical effectiveness of this approach. 

\end{abstract}

\end{frontmatter}
\textbf{Declarations of interest:} None
\section{Introduction} \label{introduction}
Signomial geometric programming (SGP) problems are optimization problems involving both positive and negative monomials in the constraints and objective function \citep{duffin1973geometric} (see \ref{sec: signomial p} for a definition). SGP problems have wide-ranging applications including engineering design \citep{avriel1978optimal,marin2007optimization, xu2013steady}), inventory control \citep{kim1998optimal,jung2005optimal,mandal2006inventory}, gas networks \citep{mishra2017comprehensive}, project management \citep{scott1995allocation}, aircraft design \citep{kirschen2018application, ozturk2019optimal}, and power control \citep{ref1corrected}.

SGP problems belong to a class of NP-hard (see \cite{xu2014global}), nonconvex, nonlinear problems, posing computational barriers to attaining global optimality (see, e.g. \citep{opgenoord2017comparison, murray2021signomial}). As such, numerous local heuristic algorithms have been proposed to solve SGP (see, e.g., \cite{kortanek1997infeasible, yang1997investigation,toscano2012heuristic,ref1corrected}). In this paper we propose a convex relaxation, suitable for either heuristic or global approaches. This places our work in the literature of convex approximations and relaxations for SGP, where a tractable convex formulation replaces the original. Such formulations are typically embedded in some sort of sequential procedure, detailed as follows.

In terms of convex approximations, a wide range of techniques have been explored. \cite{shen2005linearization} uses linear approximations to construct a sequence linear programming problems for SGP. Piecewise approximation techniques have also been used in \cite{li2009global, lundell2009a,lin2012range} to convexify signomial constraints. \cite{xu2014global} employed an inner approximation of the arithmetic-geometric (AM-GM) mean inequality and introduced an optimization algorithm through solving a sequence of geometric programming (GP) problems. \cite{li2005treating} proposed convex understimators for SGP problems with free-sign variables, utilizing exponential transformation, inverse transformation, and power convex transformation. 
Approximation via reciprocal transformation $y_j=x_j^{-\beta}$ and linear underestimation of concave terms were proposed in \cite{li2008convex} and \cite{lu2010convex} for posynomial functions. \cite{lu2012efficient} proposed a method to improve the efficiency of previous convex underestimators by finding the best value for $\beta$. However, determining an appropriate transformation function is non-trivial \citep{tseng2015milp}. In  \cite{lu2019efficient}, a convex inner approximation strategy was proposed to select an appropriate transformation.

There is also substantial literature on convex relaxations, which provide dual bounds for the (global) optimal value of SGP problems. Notably, \cite{maranas1997global} proposed an optimization algorithm that leverages convex relaxations for SGP. Linear relaxations for SGP were also proposed by \cite{shen2004global} and \cite{qu2008global}, employing exponential variable transformations and the successive refinement of these linear relaxations. Furthermore, \cite{chandrasekaran2016relative, murray2021signomial,dressler2022algebraic} introduced a hierarchy of convex relaxations for a structured class of SGP with globally nonnegative signomials, and developed non-negativity certificates for solving such SGPs.

The proposed convexification techniques, i.e., convex inner approximations or relaxations, employed for solving SGP have two limitations:

First, the effectiveness of the inner approximation approach relies on the selection of weighting factors and initial feasible solutions around which the approximations are derived \citep{hutter2011sequential}. Alternatively, since finding feasible solutions can be challenging, some variants can initialize with non-feasible solutions at the expense of additional iterations and longer time to converge \citep{xu2014global}.

Secondly, the existing convex relaxations in the literature, although independent of the choice of initial feasible solutions, are tailored for specific classes of SGP and may not be applicable to the general form of SGPs. Particularly, \cite{maranas1997global, shen2004global, qu2008global} all rely on nonzero upper and lower bounds for all variables. Further, the convex relaxations of \cite{chandrasekaran2016relative, murray2021signomial,dressler2022algebraic} are applicable only to nonnegative signomials, limiting the scope of their approach.

This paper contributes to the SGP literature and introduces novel convex relaxations for SGP, based on Exponential Conic Programming (ECP), that addresses the aforementioned limitations. ECP problems consist of a linear objective, linear constraints, and exponential cones (see \ref{sec: proposed s ecp} for details) and can be efficiently solved using interior point methods that deploy recent advancements for nonsymmetric cones \citep{nesterov2012towards, skajaa2015homogeneous, dahl2022primal, papp2022alfonso}. As of now, MOSEK \citep{aps2020mosek} is the only commercially available solver capable of effectively handling such problems. We note also that ECP can be interpreted as an generalization of GP (see, e.g., \citep{serrano2015algorithms}).

Our ECP relaxation does not require variable bounds; it can, however, be strengthened by leveraging bounds on either variables or monomials (if available) in order to add valid linear inequalities. Notably, our strengthened ECP captures the convex hull of SGP in the projection of certain variables under the settings of \cite{shen2004global,qu2008global}, where all variables possess both nonzero lower and upper bounds. Our strengthened relaxation may be viewed as analogous to the well-known SDP+RLT relaxation in polynomial optimization (see e.g. \cite{sherali2002enhancing,anstreicher2010computable}). We also embed the proposed ECP relaxation in an iterative optimization algorithm. The efficacy of the proposed relaxations (both stand-alone and within the iterative algorithm) is demonstrated through our numerical experiments.

The remainder of this paper is structured as follows: 
Section \ref{sec: signomial p} discusses SGP. In Section \ref{sec: proposed s ecp}, we introduce the ECP relaxations and valid inequalities. Section \ref{sec: geo proggramming} details sequential optimization algorithms. Numerical results are discussed in Section \ref{sec: num results}. Finally, the paper is concluded in Section \ref{sec: conclusion}.

\section{Signomial Programming Reformulation} \label{sec: signomial p}

Let $x = (x_1, \ldots, x_n) \in \mathbb{R}^n_{++}$ denote a vector of real, (strictly) positive variables. A real-valued function $f$ of $x$, with the form
\[f_k(x) = \sum_{j=1}^{m_k} c_{jk} \prod_{i=1}^n x_i^{a_{ijk}},\]
where $c_{jk} \in \mathbb{R}$ and $a_{ijk} \in \mathbb{R}$, is referred to as a polynomial function. The individual terms $g_{jk}(x) = c_{jk} \prod_{i=1}^n x_i^{a_{ijk}}$ within $f_k(x)$ are monomial functions. Furthermore, if $c_{jk} \in \mathbb{R}_+$ for all $j \in \{1, ..., n\}$ and $k \in \{1, \dots, p\}$ in $f_k(x)$, the polynomial functions $f_k(x)$ are specifically termed posynomial functions, or posynomials.

SGP can be expressed as the following non-convex nonlinear optimization problem \citep{xu2014global, kirschen2016signomial}:
\begin{subequations} \label{sgn form}
\begin{align} 
 \qquad \min \quad &  f_0(x) = \sum_{j=1}^{m_0} c_{j0} \prod_{i=1}^{n} x_i^{a_{ij0}}\\
 \text{s.t. }\quad & f_k(x)=\sum_{j=1}^{m_k} c_{jk} \prod_{i=1}^n x_i^{a_{ijk}} \leq 1, \qquad  \forall k \in \{1,\dots, p\} \label{sgp form eq1}\\
& x > 0.
\end{align}
\end{subequations}
{\color{blue}
Note that equality constraints can be reduced to \eqref{sgp form eq1}}. It is also noteworthy that in SGP \eqref{sgn form}, $c_{jk} \in \mathbb{R}$ where  some $c_{jk}$ can be negative, representing the differences among posynomials. As a special case, if we restrict $c_{jk} \in \mathbb{R}_+$ for all $j \in \{1, \ldots, n\}$ and $k \in \{1, \ldots, p\}$, SGP \eqref{sgn form} reduces to  geometric programming (GP) \citep{boyd2007tutorial}, a special case that can be solved in polynomial time.

To construct our relaxation, we develop a concise reformulation of \eqref{sgn form} by employing a series of intermediate reformulations. To this end, we first transfer the signomial objective function $f_0(x)$ from \eqref{sgn form} into the constraint set, using two strictly positive variables $x_{n+1}$ and $x_{n+2}$ as follows:
\begin{subequations} \label{sgn form inter1}
\begin{align} 
 \qquad \min \quad & x_{n+1} - x_{n+2}  \\
 \text{s.t. }\quad &
 f_0(x) = \sum_{j=1}^{m_1} c_{j0} \prod_{i=1}^{n} x_i^{a_{ij0}} \le x_{n+1} - x_{n-1} \\
 & f_k(x)=\sum_{j=1}^{m_k} c_{jk} \prod_{i=1}^n x_i^{a_{ijk}} \leq 1 \qquad  \forall k \in {1,\dots, p} \label{sgp inter1 form eq1},\\
& x_i > 0, \qquad \forall i \in  \{1, ..., n+1, n+2\}. \label{eq:xpositivity}
\end{align}
\end{subequations}


Furthermore, we split signomial terms as the difference of posynomials, in order that only posynomials are kept on either side of constraints. In particular, let $C_k^+ = \{j : c_{jk} > 0, \forall j \in M_k, k \in K \}$ and $C_k^- = \{j : c_{jk} < 0, \forall j \in M_k, k \in K \}$ be the set of monomials in constraint $k$ (including $k=0$) with positive and negative coefficients, respectively. Defining $c'_{jk}:=-c_{jk}, \forall j \in C^-$,  SGP \eqref{sgn form inter1} is reformulated as 
\begin{subequations} \label{sgn form 2}
\begin{align}   \qquad \min \quad &  x_{n+1} - x_{n+2} \\
 \text{s.t. }\quad & \sum_{j \in C_0^+} c_{j0} \prod_{i = 1 }^n x_{i}^{a_{ij0}} +   x_{n+2} \le \sum_{j \in C_0^-} c'_{j0} \prod_{i = 1 }^n x_{i}^{a_{ij0}} + x_{n+1}, \label{sgp help 1} \\ 
 & \sum_{j \in C_k^+} c_{jk} \prod_{i = 1 }^n x_{i}^{a_{ijk}}  \le \sum_{j \in C_k^-} c'_{jk} \prod_{i = 1 }^n x_{i}^{a_{ijk}} + 1  \qquad \forall k \in \{1,\dots, p\}, \label{sgp help 2} \\
& x_i > 0 \qquad \forall i \in \{1, \dots, n+2\}.
\end{align}
\end{subequations}

For simplicity in notation, we define the (extended) index set $N = \{1, 2, ..., n+1, n+2\}$ corresponding to the variable vector $x = (x_1, \dots, x_{n+2})$. Additionally, within the signomial constraint \eqref{sgp help 1}, indexed by $k=0$, we define $M_0 = \{1, \dots, m_0, m_0+1, m_0+2 \}$ as the set of monomials for this constraint. In particular, we designate $x_{n+1}$ and $x_{n+2}$ as the $(m_0+1)^{\text{th}}$ and $(m_0+2)^{\text{th}}$ monomials, respectively. These terms are then consolidated in $\sum_{j \in C_0^-} c'_{j0} \prod_{i \in N } x_{i}^{a_{ij0}}$ and $\sum_{j \in C_0^+} c_{j0} \prod_{i \in N } x_{i}^{a_{ij0}}$, respectively.
Similarly, the index set $M_k = \{1, \dots, m_k, m_k+1\}$ captures all monomials in signomial constraint \eqref{sgp help 2}, where monomial 1 is indexed by $j=(m_k+1)$ and is combined into $\sum_{j \in C_k^-} c_{jk} \prod_{i \in N } x_{i}^{a_{ijk}}$.

Given the above definitions, SGP \eqref{sgn form} can be reformulated as follows: \eqref{sgn form 3}
\begin{subequations} 
\label{sgn form 3}
\begin{align}  
\min_{x} \quad & d^T x \\
\text{s.t.} \quad & f^+_k(x) = \sum_{j \in C_k^+} c_{jk} \prod_{i \in N} x_{i}^{a_{ijk}} \le f^-_k(x) = \sum_{j \in C^-_k} c'_{jk} \prod_{i \in N} x_{i}^{a_{ijk}}, \quad \forall k \in K,  \label{sgn form 3 eq1}\\
& x > 0,
\end{align}
\end{subequations}

where $d = (0,0,...,1,-1)$ so that $d^Tx$ recovers the objective function $x_{n+1}-x_{n+2}$. Moreover, we set
\begin{equation}
\begin{cases}
c_{m_0+2,0} = 1, c'_{m_k+1,k} = 1 \qquad \qquad \quad \forall k \in K \\
a_{n+2,m_0+2,0} = 1, a_{i,m_0+2,0} = 0 \qquad \forall i \in N/\{n+2\} \\
a_{n+1,m_k+1,k} = 1 , a_{i,m_k+1,k} = 0 \qquad \forall i \in N/\{n+1\}, k \in K,
\end{cases}    
\end{equation}
to ensure equivalence between \eqref{sgn form 2} from formulation \eqref{sgn form 3}.

Formulation \eqref{sgn form 3} is a concise representation of a reformulation for SGP \eqref{sgn form}, and serves as the foundation to develop our relaxation formulation in the subsequent section.


\section{Proposed Exponential Conic Relaxations} \label{sec: proposed s ecp}
In this section, we first formally define the exponential cone. Second, we propose the ECP relaxation for SGP in general setting, requiring no bounds on nonzero variables. 
Finally, we strengthen the ECP relaxations by using bounds on either variables or monomials (if any are available) in order to add valid linear inequalities. 
\begin{definition}  \label{def. exp cone}
The exponential cone is the closure (denoted by cl) of a set of points in $\mathbb{R}^3$ given by:
\begin{equation} \label{eq: K_exp}
K_{\exp} := \mbox{cl}\{(x_1, x_2, x_3) \in \mathbb{R}^3 : x_2 \text{e}^{x_3/x_2} \le x_1, \hspace{0.1cm}x_2 \ge 0
\}.
\end{equation}
Moreover, $K_{\exp}$ can be expressed as the union of two sets $K_1 \cup K_2$, where
$K_1 := \{(x_1, x_2, x_3) \in \mathbb{R}^3 : x_2 \text{e}^{x_3/x_2} \le x_1, \hspace{0.1cm}x_2 > 0
\}$ and $K_2 := \{(x_1, 0, x_3) \in \mathbb{R}^3 : x_1\ge 0, x_3 \le 0 \}$.
\end{definition}


ECPs are defined with a linear objective over exponential cones and linear constraints. With variable transformations we can reformulate SGP \eqref{sgn form 3} such that it yields a natural ECP relaxation. To this end, let $\lambda_{jk}$ and $\gamma_{jk}$ be nonnegative auxiliary variables corresponding to the monoomials of $f_k^+(x)$ and $f_k^-(x)$, respectively. In particular, let $\lambda_{jk} \ge \prod_{i \in N } x_{i}^{a_{ijk}}, \forall j \in C^+_k, k \in K$, and $\gamma_{jk} \le \prod_{i \in N } x_{i}^{a_{ijk}}, \forall j \in C^-_k, k \in K$; moreover, introduce variables $\tilde{x}_i := \log(x_i), \forall i \in N$. We note that strict positivity of $x$ (\ref{eq:xpositivity}), as well as nonnegativity of the exponential function ensure equivalence across log transforms.  Now, observe that

\begin{equation}
\prod_{i \in N } x_{i}^{a_{ijk}} = \text{exp}\Big( \log\big(\prod_{i \in N } x_{i}^{a_{ijk}}\big) \Big) = \text{exp}\Big( \sum_{i \in N } a_{ijk} \log(x_{i}) \Big) = \text{exp}\Big( \sum_{i \in N } a_{ijk} \tilde{x}_{i} \Big). \label{change var eq1}
\end{equation}

Introducing variables $\tilde{\gamma}_{jk} := \log(\gamma_{jk})$, the inequality $\gamma_{jk} \le \prod_{i \in N } x_{i}^{a_{ijk}}$ can be reformulated as follows (applying natural logs):
\begin{equation}
\tilde{\gamma}_{jk} \le \sum_{i \in N} {a_{ijk}}\tilde{x}_{i}.\label{change var eq2} 
\end{equation}

Given \eqref{change var eq1} and \eqref{change var eq2}, an exact reformulation of SGP \eqref{sgn form 3} is obtained as follows:
\clearpage
\begin{subequations} \label{ECP SGN}
\begin{align}
\quad  \min_{x,\tilde{x}, \gamma, \tilde{\gamma}, \lambda} \quad &  d^Tx \label{eq1 ECP sgn }\\ 
\text{s.t.} \quad & \sum_{j \in C_k^+} c_{jk} \lambda_{jk} \le \sum_{j \in {C^-_k}}c'_{jk} \gamma_{jk} & \forall k \in K  \label{eq2 ECP sgn },\\
& \text{exp} \bigg(\sum_{i \in N} {a_{ijk}}\tilde{x}_{i}\bigg) \le \lambda_{jk} & \forall j \in C^+, k \in K \label{eq3 ECP sgn }, \\
& \tilde{\gamma}_{jk} \le  \sum_{i \in N} {a_{ijk}}\tilde{x}_{i}   & \forall j \in C^-, k \in K \label{eq4 ECP sgn }, \\
&\text{e}^{\tilde{x}_{i}} \le {x}_{i} \label{eq5 ECP sgn } &  \forall i \in N, \\
&{x}_{i} \le \text{e}^{\tilde{x}_{i}} \label{eq5 not-c for x} &  \forall i \in N, \\
&\text{e}^{\tilde{
\gamma}_{jk}} \le {\gamma}_{jk} &  \forall j \in C^-_k, k \in K
\label{eq6 ECP sgn }, \\
&{\gamma}_{jk} \le \text{e}^{\tilde{
\gamma}_{jk}} &  \forall j \in C^-_k, k \in K,
\label{eq6 not-c for lambda} \\
& x > 0. \label{var domain}
\end{align}
\end{subequations}
where constraints \eqref{eq5 ECP sgn } and \eqref{eq5 not-c for x}
enforce $\tilde{x}_i = \log(x_i)$; likewise, $\tilde{\gamma}_{jk} = \log(\gamma_{jk})$ is enforced in constraints \eqref{eq6 ECP sgn } and \eqref{eq6 not-c for lambda}.

In problem \eqref{ECP SGN}, constraints \eqref{eq3 ECP sgn }, \eqref{eq5 ECP sgn }, and \eqref{eq6 ECP sgn } can be represented with exponential conic constraints, following the definition of $K_{\exp}$ in \eqref{def. exp cone}. However, constraints \eqref{eq5 not-c for x} and \eqref{eq6 not-c for lambda} correspond to the hypograph of the exponential functions and are non-convex.  Relaxing these constraints yields the following natural ECP relaxation (ECPR) for SGP \eqref{sgn form 3}:
\begin{subequations} \label{ECP Relax SGP}
\begin{align}
\text{(ECPR):} \qquad \min_{x,\tilde{x}, \gamma, \tilde{\gamma}, \lambda} \quad &  d^T x  \\ 
\text{s.t.} \quad &  \eqref{eq2 ECP sgn }, \eqref{eq3 ECP sgn }, \eqref{eq4 ECP sgn }, \eqref{eq5 ECP sgn }, \eqref{eq6 ECP sgn }, \eqref{var domain}.
\end{align}
\end{subequations}
This relaxation \eqref{ECP Relax SGP} is the first convex relaxation in the literature for SGP under a general setting without explicit nonzero variable bounds. 

In the following subsections, we present valid inequalities designed to enhance the ECP relaxation~\eqref{ECP Relax SGP} with either variable or monomial bounds, provided that they are available.



\subsection{Strengthened ECP relaxation of SGP with variable upper and lower bounds} \label{subsection variable upper and lower bounds}
In this subsection, we explore a scenario in which some or all of the variables in SGP~\eqref{sgn form 3}, denoted as $x_i$ for $i \in N' = {1,...,n}$, are subject to both upper and lower bounds. Here, $N'$ represents the set of original variables in SGP \eqref{sgn form}. For a given $i \in N'$ we let $\underline{x}_i$ and $\overline{x}_i$ represent known lower and upper bounds for a variable $x_i$, such that $\underline{x}_{i} \le x_i \le \overline{x}_{i}$, and  $0 < \underline{x}_i<\overline{x}_i$. We consider two cases:

\subsubsection{Case 1: $N' \subset \{1,...,n\}$}

In this case, only a strict subset of variables have known nonzero upper and lower bounds. Consider a bounded variable, $\underline{x}_{i} \le x_{i} \le \overline{x}_{i}$ where $i \in N'$. It follows that $\log \underline{x}_{i} \le \tilde{x}_{i} \le \log \overline{x}_{i}, \forall i \in N' $. 

Let $S^x_i$ be the bounded non-convex set capturing constraint \eqref{eq5 not-c for x} and the bounds on $x_i$ and $\tilde{x}_i$ as follows
\begin{equation} \label{six}
S^x_i = \{(\tilde{x}_i, x_i): x_i \le \text{e}^{\tilde{x}_i}, \hspace{0.1cm} \log\underline{x}_{i} \le \tilde{x}_{i} \le \log\overline{{x}}_{i}, \hspace{0.1cm}  \underline{x}_{i} \le x_{i} \le \overline{x}_{i} \}, \quad \forall i \in N'.
\end{equation}
Proposition~\ref{prop cuts} establishes valid inequalities for $S^x_i$.
\begin{prop} \label{prop cuts}
For each $i \in N' \subseteq \{1,...,n\}$, the inequality
\begin{equation}
\begin{aligned}
& {x}_{i} \le \frac{\overline{x}_{i}-\underline{x}_{i}}{\log\overline{x}_{i}-\log \underline{x}_{i}}(\tilde{x}_{i}- \log \underline{x}_{i}) + \underline{x}_{i}
&    \label{ECP cut 2:initial} \qquad 
\end{aligned}
\end{equation}
is valid for $S^x_i$.
\end{prop}

\begin{proof}

By definition of $N'$, we have $ 0<\underline{x}_{i}  \le {x}_{i} \le \overline{x}_{i}$ and thus $ \log \underline{x}_{i}  \le \tilde{x}_{i} \le \log \overline{x}_{i}, \forall i \in N'$. Due to the convexity of the exponential function $\text{e}^{\tilde{x}_{i}}$, the line segment passing through points $(\log\underline{x}_{i}, \underline{x}_{i})$ and $(\log \overline{x}_{i},\overline{x}_{i})$ provides an upper bound/overestimator for $\text{e}^{\tilde{x}_{i}}$, and so
\begin{equation}
{x}_{i} \le \text{e}^{\tilde{x}_{i}} \le \frac{\overline{x}_{i}-\underline{x}_{i}}{\log\overline{x}_{i}-\log \underline{x}_{i}}(\tilde{x}_{i}- \log \underline{x}_{i}) + \underline{x}_{i}.\label{eq y2}
\end{equation}
\end{proof}


From Proposition~\ref{prop cuts} it follows that the set
\begin{equation} \label{Rx def}
R^x_i = \{(\tilde{x}_i, x_i): (3.12), \hspace{0.1cm} \tilde{x}_{i} \le \log\overline{{x}}_{i}, \hspace{0.1cm}  \underline{x}_{i} \le x_{i} \}
\end{equation}
is a convex (linear) relaxation of $S_i^x$ for all $i \in N'$. The following proposition shows that $R^x_i$ is, in fact, the convex hull of $S^x_i, \forall i \in N'$.

\begin{prop} \label{prop conv hull}
$R^x_i = \text{conv}(S^x_i) , \forall i \in N'$.
\end{prop}

\begin{proof}
For all $i \in N'$, we will show $\text{conv}(S^x_i) \subseteq R^x_i $ and $ R^x_i \subseteq\text{conv}(S^x_i) $:

Inclusion \(\text{conv}(S^x_i) \subseteq R^x_i\) follows from Proposition~\ref{prop cuts} as \(R^x_i\) is a convex relaxation of \(S^x_i\). Conversely, to show \(R^x_i \subseteq \text{conv}(S^x_i)\),  consider the polyhedron $R^x_i, \forall i \in N'$ defined by inequalities \eqref{eq1 R}-\eqref{eq4 R}:
\begin{align}
& {x}_{i} \le \frac{\overline{x}_{i}-\underline{x}_{i}}{\log\overline{x}_{i}-\log \underline{x}_{i}}(\tilde{x}_{i}- \log \underline{x}_{i}) + \underline{x}_{i}  \label{eq1 R}, \\
& \tilde{x}_{i} \le \log\overline{{x}}_{i}, \label{eq3 R} \\
& \underline{x}_{i} \le x_{i} \label{eq4 R}.
\end{align}    

We identify extreme points of \(R^x_i\) for each $i \in N'$ by setting two inequalities among \eqref{eq1 R}-\eqref{eq4 R} binding. All three bases are enumerated below:
\begin{enumerate}
\item Inequalities \eqref{eq1 R} and \eqref{eq3 R}: $(\tilde{x}_i, x_i) = (\log \overline{x}_i, \overline{x}_i)$.
\item Inequalities \eqref{eq1 R} and \eqref{eq4 R}: $(\tilde{x}_i, x_i) = (\log \underline{x}_i, \underline{x}_i)$.
\item Inequalities \eqref{eq3 R} and \eqref{eq4 R}: $(\tilde{x}_i, x_i) = (\log \overline{x}_i, \underline{x}_i)$.
\end{enumerate}

Set \(E^x_i=\{(\log \underline{x}_i, \underline{x}_i), (\log \overline{x}_i, \underline{x}_i), (\log \overline{x}_i, \overline{x}_i) \} \) captures all extreme points of \(R^x_i\) and satisfies the constraints of $S^x_i$. Thus, \(E^x_i \subseteq S^x_i\), implying \(\text{conv}(E^x_i) \subseteq \text{conv}(S^x_i)\). Given that any point \((\tilde{x}_i, x_i) \in R^x_i\) can be expressed as the convex combination of its extreme points \((E^x_i)\), it follows \(\text{conv}(E^x_i) = R \subseteq \text{conv}(S^x_i)\). Thus, it completes the proof.  
\end{proof}
\begin{figure}
\centering
\begin{tikzpicture}
\draw[->] (0,0) -- (6,0) node[right] {$\tilde{x}_i$};
\draw[->] (0,0) -- (0,6) node[above] {$x_i$};
\clip (-1,-1) rectangle (6.3,6);
\draw[blue, domain=0:6, samples=100] plot (\x, {(exp(0.35*\x)-0.2)}) node[right] {};
\draw[blue] (4.5,4.7) node[right] {$x_i = \textit{e}^{\tilde{x}_i}$};
\filldraw[black] (1.4,1.432) circle (2pt) node[above right] {};
\filldraw[black] (4.31,4.32) circle (2pt) node[above right] {};
\filldraw[black] (4.31,1.432) circle (2pt) node[above right] {};
\filldraw[red] (1.4,4.32) circle (2pt) node[above right] {};
\draw[black] (1.4,1.432) -- (4.31,4.32) node[right]{};
\draw[black] (4.31,1.432) -- (4.31,4.32) node[right]{};
\draw[black] (1.4,1.432) -- (4.31,1.432) node[right]{};
\draw[black, dashed] (1.4,1.432) -- (0,1.432) node[right]{};
\draw[black, dashed] (4.9,1.432) -- (4.31,1.432) node[right]{};
\draw[black, dashed] (4.31,0) -- (4.31,1.432) node[right]{};
\draw[black, dashed] (4.31,4.32) -- (4.31,4.9) node[right]{};
\draw[red, dashed] (0,4.32) -- (4.9, 4.32) node[right]{};
\draw[red, dashed] (1.4,0) -- (1.4,4.9) node[right]{};
\draw[black] (1.05,-0.3) node[right] {$\log \underline{x}_i$};
\draw[black] (3.95,-0.3) node[right] {$\log \overline{x}_i$};
\draw[black] (-0.5,1.432) node[right] {$ \underline{x}_i$};
\draw[black] (-0.5,4.32) node[right] {$ \overline{x}_i$};
\draw[black] (4.2,4.1) node[right] {$ (\log \overline{x}_i,\overline{x}_i)$};
\draw[black] (4.25,1.2) node[right] {$ (\log \overline{x}_i,\underline{x}_i)$};
\draw[black] (1.35,1.2) node[right] {$ (\log \underline{x}_i,\underline{x}_i)$};
\end{tikzpicture}
\caption{Representation of non-convex set $S^x_i$ and its corresponding convex hull $R^x_i$.}
\label{fig: convex hull}
\end{figure}
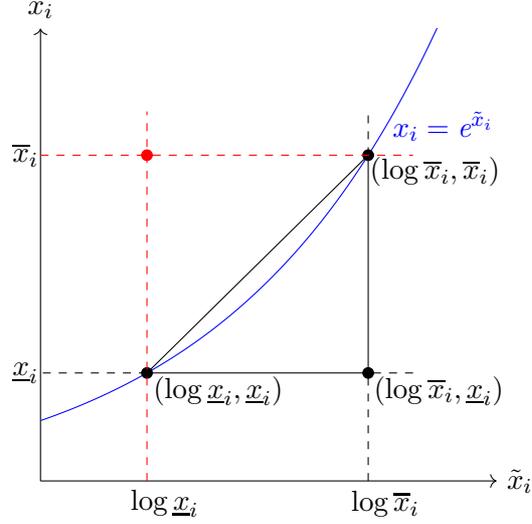

It follows from Propositions \ref{prop cuts} and \ref{prop conv hull} that the ECP relaxation \eqref{ECP Relax SGP} is strengthened by incorporating the inequalities of $R^x_{i}$. Further, Figure~\ref{fig: convex hull} represents the non-convex set $S^x_i$ alongside the corresponding convex hull formed by the inequalities of $R^x_i$. This suggests that while the bounds $\log\underline{{x}}_{i} \le \tilde{x}_{i}$ and ${x}_{i} \le \overline{{x}}_{i}$ are essential for deriving the valid inequality \eqref{eq1 R}, they are no longer explicitly required as constraints in deriving the convex hull formulation, as they become redundant to the facet-defining inequalities.




\begin{remark}
Let $\underline{x}_{i} \le x_{i} \le \overline{x}_{i}, \forall i \in N' \subset \{1,...,n\}$. The strengthened ECP relaxation (s-ECPR) of \eqref{ECP Relax SGP} is derived as follows
\begin{equation} \label{ECP S var bounds 1}
\text{s-ECPR:} \begin{cases} 
\eqref{ECP Relax SGP},  \\
(\tilde{x}_i, x_i) \in R^x_i, \qquad \forall i \in N' \in \{1,...,n\}.
\end{cases}
\end{equation}
\end{remark}

\subsubsection{Case 2: $N'=\{1,...,n\}$}
In this case, all original variables are bounded from below and above, characterized by $\underline{x}_{i} \le x_{i} \le \overline{x}_{i}, \forall i \in N'=\{1, \dots, n\}$, where $0 < \underline{x}_{i}< \overline{x}_{i}$. To enhance the ECP relaxations \eqref{ECP Relax SGP} and \eqref{ECP S var bounds 1}, we derive implicit bounds for the auxiliary variables $x_{n+1}$, $x_{n+2}$, and $\gamma$.


By construction from \eqref{sgn form 2}, it is clear that $x_{n+1}-x_{n+2}$ capture the objective function value $\sum_{j \in C^+_0} c_{j0} \prod_{i=1}^{n} x_i^{a_{ij0}} - \sum_{j \in C^-_0} c'_{j0} \prod_{i=1}^{n} x_i^{a_{ij0}}$. Similarly, as stated in previous section, $\gamma_{jk}$ corresponds to the value of $\prod_{i \in N}{x}^{a_{ijk}}$. For all $j \in C_k^-, k \in K$, letting $A^+(j) = \{i : a_i \ge 0, i=1,\dots, n \}$ and $A^-(j) = \{i : a_i < 0, i ={1,\dots, n} \}$, we can easily establish nonzero lower and upper bounds for $x_{n+1}$, $x_{n+2}$, and $\gamma$ as follows
\begin{equation} \label{additional bounds}    
\begin{cases}
& \underline{x}_{n+1} = \sum_{j \in C^+_0} c_{j0} \prod_{i\in A^+(j)}\underline{x}_i^{a_{ijk}}\prod_{i\in A^-(j)}\overline{x}_i^{a_{ijk}}, \\
& \overline{x}_{n+1} = \sum_{j \in C^+_0} c_{j0} \prod_{i\in A^+(j)}\overline{x}_i^{a_{ijk}}\prod_{i\in A^-(j)}\underline{x}_i^{a_{ijk}}, \\
& \underline{x}_{n+2} =  \sum_{j \in C^-_0} c_{j0} \prod_{i\in A^+(j)}\underline{x}_i^{a_{ijk}}\prod_{i\in A^-(j)}\overline{x}_i^{a_{ijk}}, \\
& \overline{x}_{n+2} =  \sum_{j \in C^-_0} c_{j0} \prod_{i\in A^+(j)}\overline{x}_i^{a_{ijk}}\prod_{i\in A^-(j)}\underline{x}_i^{a_{ijk}}, \\
& \underline{{\gamma}}_{jk} = \prod_{i\in A^+(j)}\underline{x}_i^{a_{ijk}}\prod_{i\in A^-(j)}\overline{x}_i^{a_{ijk}} \qquad \forall j \in C_k^-, k \in K,  \\
& \overline{{\gamma}}_{jk} = \prod_{i\in A^+(j)}\overline{x}_i^{a_{ijk}}\prod_{i\in A^-(j)}\underline{x}_i^{a_{ijk}} \qquad \forall j \in C_k^-, k \in K.
\end{cases}
\end{equation}

Using the established bounds for variables $x_{n+1}$ and $x_{n+2}$, we can derive the valid inequalities of $R^x_i$ in \eqref{Rx def} for these variables, thereby improving the relaxation associated with the constraints \eqref{eq5 not-c for x}.

Further, the bounded non-convex set corresponding to constraint \eqref{eq6 not-c for lambda} is expressed as
\begin{equation}
\begin{aligned}
 S^{\gamma}_{jk} = \Big\{(\tilde{\gamma}_{jk}, \gamma_{jk}): & \quad \gamma_{jk} \le \text{e}^{\tilde{\gamma}_{jk}}, \\ 
& \quad \underline{{\gamma}}_{jk} \le {\gamma}_{jk} \le \overline{{\gamma}}_{jk}, \\
 &
\hspace{0.1cm}\log(\underline{{\gamma}}_{jk}) \le {\tilde{\gamma}}_{jk} \le \log (\overline{{\gamma}}_{jk}) \Big\}, \quad  \forall j \in C_k^-, k \in K.     
\end{aligned}
\end{equation}

Proposition~\ref{prop conv hull lambda} finds the convex hull of $S^{\gamma}_{jk}, \forall j \in C_k^-, k \in K $
\begin{prop} \label{prop conv hull lambda}
Set 
\begin{equation} \label{R gamma def}
\begin{aligned}
 R^{\gamma}_{jk} = \Big\{(\tilde{\gamma}_i, \gamma_i) : \hspace{0.1cm} & {\gamma}_{jk} \le \frac{\overline{\gamma}_{jk}-\underline{\gamma}_{jk}}{\log\overline{\gamma}_{jk}-\log \underline{\gamma}_{jk}}(\tilde{\gamma}_{jk}- \log \underline{\gamma}_{jk}) + \underline{\gamma}_{jk}, \\
 &
{\tilde{\gamma}}_{jk} \le \log \overline{\gamma}_{j,k}, \hspace{0.1cm} \underline{\gamma}_{j,k} \le {\gamma}_{jk} \hspace{0.1cm} \Big\}  \end{aligned}
\end{equation}
is the convex hull of $S^{\gamma}_{jk},   \forall j \in C_k^-, k \in K  $. 
\end{prop}

\begin{proof}
First, it follows from Proposition~\ref{prop cuts} that $ R^{\gamma}_{jk}$ is a convex relaxation of $S^{\gamma}_{jk}, \forall j \in C_k^-, k \in K$. Further, considering the inequalities defining polyhedron $R^{\gamma}_{jk}$, set $$E^{\gamma}_{jk}=\{(\log \underline{\gamma}_{j,k}, \underline{\gamma}_{j,k}), (\log \overline{\gamma}_{j,k}, \underline{\gamma}_{j,k}), (\log \overline{\gamma}_{j,k}, \overline{\gamma}_{j,k})\}$$ contains all the extreme points of $R^{\gamma}_{jk}$. Since  $E^{\gamma}_{jk} \subseteq S^{\gamma}_{jk}$, similar to the proof of Proposition~\ref{prop conv hull} it follows that  $R^{\gamma}_{jk}=\text{conv}(S^{\gamma}_{jk}), \forall j \in C_k^-, k \in K.$   
\end{proof}

\begin{remark} \label{remark conv}
Let $\underline{x}_{i} \le x_{i} \le \overline{x}_{i}, \forall i \in N' = \{1,...,n\}$. The strengthened ECPR (s-ECPR) \eqref{s-ECP Relax SGP} includes the convex hull representation of problem \eqref{ECP SGN}.
\begin{equation} \label{s-ECP Relax SGP}
\text{s-ECPR:}    
\begin{cases}
\eqref{ECP Relax SGP}  \\
(\tilde{x}_i, x_i) \in R^x_i, \quad \qquad \forall i \in N, \\
(\tilde{\gamma}_{jk},\gamma_{jk}) \in R^{\gamma}_{jk}, \qquad \forall j \in C_k^-, k \in K.
\end{cases}    
\end{equation}

\end{remark}



\subsection{Strengthened ECP relaxation with nonzero monomial bounds} \label{subsection variable lower bounds}
 In this section, we address scenarios where monomial bounds can be established from the bounds of variables and present valid inequalities for such settings. Namely, we consider the case where some monomial terms $g_{jk}(x)=c_{jk} \prod_{i \in N } x_{i}^{a_{ijk}}$ 
are bounded either from above or below and utilize the structure of the problem in two certain cases as follows:
\subsubsection{Case 1: Monomial lower bounds}
In signomial constraint set \eqref{sgp form eq1}, let $K' \subseteq \{1, \ldots, p\}$ be the set of monomial constraints $g_{1k}(x) = c_{1k} \prod_{i \in N} x_{i}^{a_{i1k}} \le 1$ for which there exists a lower bound, represented by $\underline{X}_{k}$. Index $j=1$ indicates that there is only one monomial term in such constraints. Therefore, the monomial constraints can be written as 
\begin{equation} \label{eq monomials}
\underline{X}_{k} \le g_{1k}(x) = c_{1k} \prod_{i \in N } x_{i}^{a_{i1k}} \le 1, \forall k \in K',    
\end{equation}
where such monomial lower bounds exist only if there are lower bounds for variables with nonzero $a_{i1k}$ for all $k \in K' \subseteq K$ and are obtained as $\underline{X}_{k} = c_{1k} \prod_{i \in N} \underline{x}_{i}^{a_{i1k}}$ for all $k \in K'$.

We now reformulate monomial constraints to derive valid inequilites. Introducing non-negative variables $X_{k}$ and $\tilde{X}_{k}$, let $X_{k} = c_{1k} \prod_{i \in N } x_{i}^{a_{i1k}}$ and $\tilde{X}_{1k} = \log X_{k}$. An extended formulation of monomial constraints \eqref{eq monomials} is given by:
\begin{subequations} \label{extended monomial}
\begin{align}
& \underline{X}_{k} \le X_{k} \le 1, & \forall k \in K'  \label{monomial ref eq 1}, \\
& \text{exp}\big(  \tilde{X}_k \big) \le X_k , & \forall k \in K' \label{monomial ref eq 2},\\
& X_k \le \text{exp}\big(  \tilde{X}_k\big), & \forall k \in K' \label{monomial ref eq 3}, \\
& \tilde{X}_k =  c_{1k}\sum_{i \in N} {a_{i1k}}\tilde{x}_{i}, & \forall k \in K'  \label{monomial ref eq 4}, \\
& \log (\underline{X}_{k}) \le \tilde{X}_{k} \le 0, & \forall k \in K'. \label{bounds big X tilde}
\end{align}
\end{subequations}
Constraints \eqref{monomial ref eq 2}-\eqref{monomial ref eq 4} enforce $X_k = \text{exp}(\tilde{X}_k) = c_{1k} \prod_{i \in N } x_{i}^{a_{i1k}}$, and \eqref{bounds big X tilde} is implied from \eqref{monomial ref eq 1} and the definition of $\tilde{X}_{k}$. It follows from Proposition~\ref{prop cuts} that we can derive valid inequalities for the non-convex constraint \eqref{monomial ref eq 3} with bounded variables $X_k$ and $\tilde{X}_k$ as follows:
\begin{equation} \label{aggregated cut}
{X}_{k} \le \frac{\underline{X}_{k}-1}{\log \underline{X}_{k}}\bigg(\tilde{X}_{k} - \log \underline{X}_{k}\bigg) + \underline{X}_{k}, \qquad \forall k \in K'.
\end{equation}
\begin{remark}
The strengthened SGP relaxation with monomial lower bounds is obtained as follows:
\begin{equation}\label{s-ECPR 3}
\text{s-ECPR:}
\begin{cases}
\eqref{ECP Relax SGP},  \\
\eqref{monomial ref eq 1}, \eqref{monomial ref eq 2}, \eqref{monomial ref eq 4}, \eqref{bounds big X tilde},\eqref{aggregated cut}.
\end{cases}    
\end{equation}    
\end{remark}

\subsubsection{Case 2. Monomial upper bounds} \label{subsection variable upper bounds}
In \eqref{sgn form 3 eq1}, we have $f^-_k(x) =  \sum_{j \in {C^-_k}} g_{jk}$ where $g_{jk}=c'_{jk} \prod_{i \in N } x_{i}^{a_{ijk}}$ are monomial terms. Let $\overline{X}_{jk} = c'_{jk}\prod_{i \in N} \overline{x}_{i}^{a_{ijk}}$ represent the upper bound of monomial terms $g_{jk}(x)$ for $j$ in a subset of $C^-_k$, denoted as $J_k$, and $k \in K$. Such upper bound exists, for a given $j \in C_k^-, k \in K' \subseteq K$, 
only if there exist $\overline{x}_t$ for all $t \in T_{jk} = \{i\in N : a_{ijk} \neq 0 \}$. The following proposition describes valid inequalities for SGP of \eqref{ECP SGN} using these upper bounds.
\begin{prop}
Let $\overline{X}_{jk}= c'_{jk} \prod_{i \in N} \overline{x}_{i}^{a_{ijk}}$ be the upper bound of monomial $g_{jk}(x)$ for $j \in J_k \subseteq C^-_k$ and $k \in K' \subset K$. For a given $k \in K'$,
inequality 
\begin{equation}
\gamma_{jk} \le \frac{\overline{X}_{jk}}{c'_{jk}}, \qquad \forall j \in J_k \subseteq C^-_k, k \in  K' \subset K 
\label{valid ieq 1 mu}
\end{equation}
is valid for SGP \eqref{ECP SGN}. Further, if $J_k = C^-_k$, the following inequality 
\begin{equation}
\sum_{j \in {C^+_k}} c_{jk}\lambda_{jk} \le \sum_{j \in {C^-_k}}  \overline{X}_{jk}, \qquad \forall k \in K' \subset K  \label{valid ieq 2 mu}
\end{equation}
is also valid for SGP \eqref{ECP SGN}.
\end{prop}

\begin{proof}
Given the definition of $\gamma_{jk}$ in \eqref{ECP SGN}, we have
\[ \gamma_{jk} \le \prod_{i \in N } x_{i}^{a_{ijk}} \le \prod_{i \in N } \overline{x}_{i}^{a_{ijk}} = \frac{\overline{X}_{jk}}{c'_{jk}}, \qquad \forall j \in J_k \subseteq C^-_k, k \in  K' \subset K. \]

This implies that inequality \eqref{valid ieq 1 mu} is valid for SGP \eqref{ECP SGN}. Further, from \eqref{eq2 ECP sgn }, it follows that if $J_k = C^-_k$ for a $k \in K' \subset K$, we have
\[ \sum_{j \in {C^+_k}} c_{jk}\lambda_{jk} \le \sum_{j \in {C^-_k}} c'_{jk}\gamma_{jk} \le \sum_{j \in {C^-_k}} \overline{X}_{jk}.  \]

This demonstrates the validity of inequality \eqref{valid ieq 2 mu}.
\end{proof}
\begin{remark}
The strengthened SGP relaxation with monomial upper bounds is obtained as follows:
\begin{equation} \label{s-ECPR 4}
\text{s-ECPR:}
\begin{cases}
\eqref{ECP Relax SGP},  \\
\eqref{valid ieq 2 mu}, \eqref{valid ieq 1 mu}.
\end{cases}    
\end{equation}    
\end{remark}

\section{Sequential Optimization Algorithms for SGP} \label{sec: geo proggramming}

The sequential optimization algorithms in the literature generally aim to solve SGP \eqref{sgn form 3} by solving a sequence of subproblems where the non-convex terms in $f^-_k(x)$ are convexified. 
Given an initial solution $x^{(0)}$, SGP is iteratively solved through a sequence of subproblems at each iteration $t$ as follows:
\begin{subequations} \label{GP sub pro}
\begin{align}  
\text{(sub $t$):} \qquad \min_{x} \quad &  d^T x \\
 \text{s.t. }\quad & f^+_k(x) \le \hat{f}^-_k({x;x^{(t-1)}}),\qquad \forall k \in K,  \\
& x > 0.
\end{align}
\end{subequations}

In a given subproblem $t$, we have the solution from the previous iteration $x^{(t-1)}$ from which we derive $\hat{f}^-_k({x;x^{(t-1)}})$ in order to convexify $f^-_k(x)$ around the point $x^{(t-1)}$. The construction of the convex term $\hat{f}^-_k({x;x^{(t-1)}})$ varies across different methodologies, primarily relying on inner approximations. Notably, approaches such as \cite{li2009global, lundell2009a, lin2012range} employ piece-wise inner approximations for $\hat{f}^-_k({x;x^{(t-1)}})$. However, in the work of \citep{xu2014global}, $\hat{f}^-_k({x;x^{(t-1)}})$ is derived as the local monomial estimator of the posynomial ${f}^-_k({x})$ at the point $x^{(t-1)}$, utilizing Taylor's approximation. Subsequently, these algorithms proceed by iteratively searching for new solutions, with termination occurring when $| x^{(t)} - x^{(t-1)} | \le \epsilon$, where $\epsilon > 0$ denotes a specified solution accuracy.

Although \cite{xu2014global}'s sequential optimization algorithm is effective, its reliance on an arbitrary initial solution $x^{(0)}$ can potentially result in increased iterations and increased solving times, especially when such a solution is infeasible or of poor quality \citep{xu2014global}. In contrast, we present adapt \cite{xu2014global}'s framework to our ECP relaxation and avoid such issues. At iteration $t=0$, our proposed relaxation is used for an initial solution, denoted as $\tilde{x}^{(0)}$ and $\tilde{\gamma}^{(0)}$. Subsequently, for each iteration $t \ge 1$, we solve an ECP subproblem, represented as \eqref{ECP sub SGN}, incorporating our proposed valid inequalities to enhance the subproblem formulation, leveraging the solutions obtained from the previous iteration.

\subsection{Proposed Sequential Optimization Algorithm}
Let $\tilde{x}^{(0)}$ and $\tilde{\gamma}^{(0)}$ be the initial solution obtained from the proposed relaxations. Formulation \eqref{ECP sub SGN} is 
subproblem $t$ given the optimal solution at iteration $t-1$, i.e., $\tilde{x}^{(t-1)}$, $\tilde{\gamma}^{(t-1)}$:
\begin{subequations} \label{ECP sub SGN}
\begin{align}
\text{(ECP sub-$t$):} \quad \min_{x,\tilde{x}, \gamma, \tilde{\gamma}, \lambda, \eta, \eta'} \quad &  d^T x  + w \eta + w' \eta'\label{eq1 sub ECP sgn }\\ 
 \text{s.t.} \quad & \eqref{eq2 ECP sgn }, \eqref{eq3 ECP sgn }, \eqref{eq4 ECP sgn }, \eqref{eq5 ECP sgn }, \eqref{eq6 ECP sgn }, \\
& \eqref{Rx def}, \eqref{R gamma def}, \\
& {x}_{i} \le \hat{f}(\tilde{x};\tilde{x}^{(t-1)}) + \eta_i , \quad \forall i \in N, \label{eq6 sub affine x} \\
& \gamma_{jk} \le \hat{f}(\tilde{\gamma};\tilde{\gamma}^{(t-1)}) + \eta'_{jk}, \qquad \forall j \in M_k, k \in K,
\label{eq6 sub affine gamma} \\
& \eta, \eta' \ge 0. \label{eta domain}
\end{align}
\end{subequations}
In a subproblem \eqref{ECP sub SGN}, we incorporate the inequalities from \eqref{Rx def} and \eqref{R gamma def}, assuming that all variables have upper and lower bounds. In scenarios where this assumption does not hold, these two constraints can be substituted with valid inequalities derived for alternative cases.

Further, constraints \eqref{eq6 sub affine x} and \eqref{eq6 sub affine gamma} convexify non-convex constraints \eqref{eq5 not-c for x} and \eqref{eq6 not-c for lambda}, using affine estimators 
\begin{align}
& \hat{f}(\tilde{x};\tilde{x}^{(t-1)}) = \text{e}^{\tilde{x}^{(t-1)}} + \text{e}^{\tilde{x}^{(t-1)}}(\tilde{x}-\tilde{x}^{(t-1)}),    \\
& \hat{f}(\tilde{\gamma};\tilde{\gamma}^{(t-1)}) = \text{e}^{\tilde{\gamma}^{(t-1)}} + \text{e}^{\tilde{\gamma}^{(t-1)}} (\tilde{\gamma}-\tilde{\gamma}^{(t-1)})
\end{align}
to approximate the functions $f(\tilde{x})=\text{e}^{\tilde{x}}$ and $f(\tilde{\gamma})=\text{e}^{\tilde{\gamma}}$ centered at points $\tilde{x}^{(t-1)}$ and $\tilde{\gamma}^{(t-1)}$, respectively. 

The Taylor's approximation implies $\hat{f}$ is an underestimator function so that  $x \le \hat{f}(\tilde{x};\tilde{x}^{(t-1)}) \le f(\tilde{x})$ and $\gamma \le \hat{f}(\tilde{\gamma};\tilde{\gamma}^{(t-1)}) \le f(\tilde{\gamma})$. Furthermore, constraints \eqref{eq5 ECP sgn } and \eqref{eq6 ECP sgn } enforce $f(\tilde{x}) \le x$ and $f(\tilde{\gamma}) \le \gamma$. Thus, to ensure the feasibility of subproblem \eqref{ECP sub SGN} given these underestimators, auxiliary variables $\eta \ge 0$ and $\eta' \ge 0$ are introduced in constraints \eqref{eq6 sub affine x} and \eqref{eq6 sub affine gamma}. These variables are integrated into the objective function with penalty weight matrices $\omega$ and $\omega'$, respectively. It is noteworthy that, with the penalty terms, these auxiliary variables converge to 0 at optimality.

Algorithm 1 presents the steps of the proposed iterative algorithm.
\begin{algorithm}
\caption{Sequential ECP Algorithm}
\begin{enumerate}
    \item \textbf{Initial Relaxation:}
    \begin{enumerate}
        \item Set iteration counter $t=0$.
        \item Obtain $\tilde{x}^{(0)}$ and $\tilde{\gamma}^{(0)}$ through solving ECPR \eqref{ECP Relax SGP} or s-ECPR \eqref{ECP S var bounds 1}, \eqref{s-ECP Relax SGP}, \eqref{s-ECPR 3}, or \eqref{s-ECPR 4} depending on problem structure.
        \item Set solution accuracy $\epsilon$.
    \end{enumerate}

    \item \textbf{Iterative Subproblems:}
    \begin{enumerate}
       \item  Update $t \leftarrow t+1$. 
        \item Solve ECP sub-$t$ problem \eqref{ECP sub SGN} and obtain $\tilde{x}^{(t)}$ and $\tilde{\gamma}^{(t)}$.
    \end{enumerate}
    \item \textbf{Convergence Check}
    \begin{enumerate}
    \item If $\| \tilde{x}^{(t)} - \tilde{x}^{(t-1)}\| \le \epsilon $ and $\|\tilde{\gamma}^{(t)} - \tilde{\gamma}^{(t-1)}\| \le \epsilon $, stop. 
    \item Else, go to step 2.(a).
    \end{enumerate}
\end{enumerate}
\end{algorithm}

Using the first-order Taylor approximation, \citep{boyd2007tutorial} demonstrated that finding a monomial approximation of $f^-_k(x)$ in the vicinity of $x^{(t-1)}$ in \cite{xu2014global}'s algorithm (i.e., subproblem \eqref{GP sub pro}) is equivalent to searching for an affine approximation of $\text{e}^{\tilde{x}}$ and $\text{e}^{\tilde{\gamma}}$ near $\tilde{x}^{(t-1)}$ and $\tilde{\gamma}^{(t-1)}$, respectively, in our proposed subproblem \eqref{ECP sub SGN}. Thus the iterative algorithm we propose retains the convergence properties of \cite{xu2014global}'s algorithm to a KKT point. Moreover, through numerical validation in the next section, we demonstrate that fewer iterations are needed thanks to our tightened relaxations.

\section{Numerical Results} \label{sec: num results}
In this section, we demonstrate the effectiveness of our proposed relaxations in comparison to existing methods. Subsequently, we apply the iterative algorithm to solve a set of test SGP instances. All instances are implemented in version 1.17.0 of the Julia modeling programming language (JuMP) \cite{kwon2019julia} and solved using MOSEK \cite{aps2020mosek} version 10.1.3 on a 1.5 GHz laptop with 16 GB RAM. MOSEK is chosen as the solver due to its status as the sole globally convergent software available for ECP problems.

\subsection{Evaluation of Proposed Relaxations}
In this section, we compare the quality of our proposed relaxation method by the existing relaxations on benchmark test instances used in \cite{qu2008global}. 

We provide an overview of the test sets and examine the resulting relaxations in detail for Problem 1. Subsequently, complete results are presented for other instances.

\noindent \textbf{Problem 1} \citep{qu2008global}:
\begin{align*}
\text{(P1)} \quad \min \quad & 6x_1^2 + 4x_2^2 - 2.5x_1x_2 \\
\text{s.t} \quad & -x_1x_2 \le -8, \\
&  1 \le x_1 \le 10, \hspace{0.1cm} 1 \le x_2 \le 10.
\end{align*}

\noindent \textbf{Problem 2} \citep{rijckaert1978comparison}
\begin{align*}
\text{(P2)} \quad \min \quad & 168x_1x_2 + 3651.2x_1 x_2 x_3^{-1} + 40000x_4^{-1} \\
\text{s.t.} \quad & 1.0425x_1x_2^{-1} \leq 1, \\
& 0.00035x_1x_2 \leq 1, \\
& 1.25x_1^{-1}x_4 + 41.63x_1^{-1} \leq 1, \\
& 40 \le x_1 \le 44, \hspace{0.1cm} 40 \le x_2 \le 45, \\
& 60 \le x_3 \le 70, \hspace{0.1cm} 0.1 \le x_4 \le 1.4.
\end{align*}

\noindent \textbf{Problem 3} \citep{floudas2013handbook}
\begin{align*}
\text{(P3)} \quad \min \quad & 0.4x_1^{0.67}x_7^{-0.67} + 0.4x_2^{0.67}x_8^{-0.67} + 10 - x_1 - x_2 \\
\text{s.t.} \quad & 0.0588x_5x_7 + 0.1x_1 \leq 1, \\
& 0.0588x_6x_8 + 0.1x_1 + 0.1x_2 \leq 1, \\
& 4x_3x_5^{-1} + 2x_3^{-0.71}x_5^{-1} + 0.0588x_3^{-1.3}x_7 \leq 1, \\
& 4x_4x_6^{-1} + 2x_4^{-0.71}x_6^{-1} + 0.0588x_4^{-1.3}x_8 \leq 1, \\
& 0.1 \leq x_i \leq 10, \quad i = 1, 2, \ldots, 8.
\end{align*}

\noindent \textbf{Problem 4} \citep{floudas2013handbook}

\begin{align*}
\text{(P4)} \quad \min \quad & x_1+x_2+x_3 \\
\text{s.t.} \quad & 833.33252x_1^{-1}x_4x_6^{-1} -833.333x_1^{-1}x_6^{-1} + 100 x_6^{-1}  \leq 1, \\
& 1250(x_2^{-1}x_5x_7^{-1}-x_2^{-1}x_4x_7^{-1})+x_4x_7^{-1} \leq 1, \\
& 1250000x_3^{-1}x_8^{-1} -2500x_3^{-1}x_5x_8^{-1} + x_5x_8^{-1}  \leq 1, \\
& 0.0025(x_4 + x_6) \leq 1, \\
& 0.0025(-x_4 + x_5+x_7) \leq 1, \\
& 0.01(x_8 - x_5) \leq 1, \\
& 100 \leq x_1 \leq 10000, \\
& 1000 \leq x_i \leq 10000 \quad i = 2,3, \\
& 10 \leq x_i \leq 1000 \quad i = 4,5, \dots, 8. 
\end{align*}

\noindent \textbf{Problem 5} \citep{qu2008global}:
\begin{align*}
\text{(P5)} \quad \min \quad & 5x_1 + 50000x_1^{-1} + 46.2x_2 + 72000x_2^{-1} + 144000x^{-1}_3, \\
\text{s.t.} \quad & 4x_1^{-1} + 32x_2^{-1} + 120x_3^{-1} \leq 1, \\
& 1 \leq x_1, \, x_2, \, x_3 \leq 220.
\end{align*}

\noindent \textbf{Problem 6} \citep{qu2008global}:
\begin{align*}
\text{(P6)} \quad \min \quad  & 5.3578x_3^2 + 0.8357x_1x_5 + 37.2392x_1, \\
\text{s.t.} \quad & 0.00002584x_3x_5 - 0.00006663x_2x_5 - 0.0000734x_1x_4  \le 1, \\
 & 0.00085307x_2x_5 + 0.00009395x_1x_4 - 0.00033085x_3x_5 \le 1, \\
 & 1330.3294x^{-1}_2x^{-1}_5 - 0.42x_1x^{-1}_5 - 0.30586x^{-1}_2 x^2_3 x^{-1}_5 \le 1, \\
& 0.00024186x_2x_5 + 0.00010159x_1x_2 + 0.00007379x_3^2 \le 1, \\
 & 2275.1327x^{-1}_3x^{-1}_5 - 0.2668x_1x^{-1}_5 - 0.40584x_4x^{-1}_5 \le 1, \\
& 0.00029955x_3x_5 + 0.00007992x_1x_3 - 0.00012157x_3x_4 \le 1, \\
 &78.0 \leq x_1 \leq 102.0, \hspace{0.1cm} 33.0 \leq x_2 \leq 45.0, \hspace{0.1cm} 27.0 \leq x_3 \leq 45.0, \\
&27.0 \leq x_4 \leq 45.0, \hspace{0.1cm} 27.0 \leq x_5 \leq 45.0.
\end{align*}

\noindent \textbf{Problem 7} \citep{qu2008global}:
\begin{align*}
\text{(P7)} \quad \min \quad & 0.5x_1x_2^{-1} - x_1 - 5x_2^{-1} \\
\text{s.t.} \quad & 0.01x_2x_3^{-1} + 0.01x_2 + 0.0005x_1x_3 \leq 1, \\
& 70 \leq x_1 \leq 150, \hspace{0.1cm} 1 \leq x_2 \leq 30, \hspace{0.1cm} 0.5 \leq x_3 \leq 21.    
\end{align*}

In problem 1, the optimal solution is $x^*= (2.5558, 3.1302)$, which yields an objective function value of $z^*=58.38488$. Letting $\tilde{x}_i = \log(x_i) \forall i=1,2$, Table \ref{tab: relaxations} presents three relaxations produced by the methods of \cite{maranas1997global}, \cite{shen2004global}, and \cite{qu2008global} for problem P1. The relaxation of \cite{maranas1997global} can be solved as an ECP in MOSEK by transferring the exponential terms of the objective function into constraints. In Table \ref{tab: relaxations}, $LB_i$ represents the optimal solution of each relaxation, providing a lower bound of the optimal objective function value of problem P1. As observed from the table, the linear relaxation of \cite{qu2008global} provides a tighter lower bound, i.e., $LB_3$, for problem P1.

\begin{table*}[t!]
\caption{Detained formulation of various relaxations for problem P1}
\label{tab: relaxations}
\scalebox{0.85}{
\begin{tabular}{|l|l|l|}
    \hline
   Relaxation of \cite{qu2008global}  &    Relaxation of \cite{shen2004global} &    Relaxation of \cite{maranas1997global}  \\ \hline
$\min \hspace{0.1cm}  z_3 = 1.0235\tilde{x}_1 + 5.9967\tilde{x}_2 $ &  $\min \hspace{0.1cm} z_2= 204.227\tilde{x}_1 + 118.237\tilde{x}_2$ & $\min \hspace{0.1cm} z_1=6 \text{e}^{2\tilde{x}_1} + 4\text{e}^{2\tilde{x}_2}+5$   \\ 
$\hspace{0.6cm} + 40.997$ & $\hspace{0.6cm} - 28.290$ & $\hspace{0.6cm} - 53.744(\tilde{x}_1+\tilde{x}_2)$ \\
$\text{s.t:} \hspace{0.1cm}-1.001\tilde{x}_1 - 2.995\tilde{x}_2 \le -6.0172$ & $ \text{s.t:} \hspace{0.1cm}-21.4976(\tilde{x}_1+\tilde{x}_2)\le -7$ & $\text{s.t:} \hspace{0.1cm} -21.4976(\tilde{x}_1+\tilde{x}_2)\le -6$ \\
$\hspace{0.7cm} 0 \le \tilde{x}_1 \le \log(10)$ & $\hspace{0.7cm} 0 \le \tilde{x}_1 \le \log(10)$& $\hspace{0.7cm} 0 \le \tilde{x}_1 \le \log(10)$ \\
$\hspace{0.7cm} 0 \le \tilde{x}_2 \le \log(10)$ & $\hspace{0.7cm} 0 \le \tilde{x}_2 \le \log(10)$ & $\hspace{0.7cm} 0 \le \tilde{x}_2 \le \log(10)$ \\ \hline
$\hspace{0.7cm}\tilde{x}^*=(2.30258, 1.23943)$ & $\hspace{0.7cm}\tilde{x}^*=(0.0, 0.32562)$& $\hspace{0.7cm}\tilde{x}^*=(0.9,1)$  \\ \hline 
$\hspace{0.7cm} LB_3 = 50.7862$ & $\hspace{0.7cm} LB_2 = 10.2094 $ & $\hspace{0.7cm} LB_1 =  -31.2595$ \\ \hline
\end{tabular}
}
\end{table*}

We also derive our ECP relaxation of P1 based on \eqref{ECP Relax SGP}. First, given the variable domain $x_1$ and $x_2$, the objective function of P1 is non-negative, requiring only one auxiliary variable $x_3 > 0$ to transfer the term into the constraints. Introducing variable vectors $\lambda, \gamma, \tilde{\gamma}$ with appropriate dimensions, we derive the proposed ECP relaxation for problem P1 as follows:
\begin{subequations} \label{p0 ecp}
\begin{align}
(\text{ECPR}): \quad \min \quad & x_3 \\   
\text{s.t.} \quad & 6 \lambda_{11} + 4\lambda_{12} \le 2.5\gamma_{13} + \gamma_{14},   \\ 
& 8 \le \gamma_{21}, \\
& \text{e}^{2\tilde{x}_1} \le \lambda_{11}, \\
& \text{e}^{2\tilde{x}_2} \le \lambda_{12}, \\
& \tilde{\gamma}_{13} \le \tilde{x}_1 + \tilde{x}_2, \\
& \tilde{\gamma}_{14} \le \tilde{x}_3, \\
&\tilde{\gamma}_{21} \le \tilde{x}_1+ \tilde{x}_2, \\
& \text{e}^{\tilde{x}_i} \le x_i \qquad \forall i =1,2,3, \\
& \text{e}^{\tilde{\gamma}_{1j}} \le \gamma_{1j} \qquad \forall j=3,4, \\
& \text{e}^{\tilde{\gamma}_{2j}} \le \gamma_{2j} \qquad \forall j=1, \\
& 1 \le x_i \le 10 \quad \forall i=1,2, \\
& 0 \le \tilde{x}_i \le \log(10) \quad \forall i=1,2.
\end{align}
\end{subequations}
The ECPR problem \eqref{p0 ecp} yields the optimal solution $x^*=(2.4817, 1.000, 7.4998)$, $\tilde{x}^*=(0.3791, 0, 2.0149)$, $\lambda^*=(2.4783, 1.4279)$, $\gamma^*=(5.2324, 7.5, 100)$, and $\tilde{\gamma}^*=(0.2942, 2.0149, 0.3791)$, providing the lower bound value $LB_4=z^*=7.4998$ for the objective function of P1.

Using the constraints derived from \eqref{additional bounds}, we establish bounds for auxiliary variables: $x_3 \in [7.5, 750]$, $\gamma_{13} \in [1,100]$, $\gamma_{14} \in [7.5,750]$, and $\gamma_{21} \in [1,100]$. Subsequently, we construct s-ECPR \eqref{p0 ecp} for P1 by incorporating valid inequalities \eqref{Rx def} and \eqref{R gamma def} given the above boudns. The s-ECPR relaxation of problem P1 yields an optimal solution of $x^*=(2.62, 3.18, 56.7598)$, $\tilde{x}^*=(0.5018, 0.6376, 3.3943)$, $\lambda^*=(3.0987, 4.0934)$, $\gamma^*=(9.7002, 10.7157, 100)$, and $\tilde{\gamma}^*=(0.79327, 2.1049, 1.1394)$. This results in a substantially improved lower bound value of $LB_5=z^*=56.7598$ for the objective function of P1. 

\begin{table}
\caption{Comparison of various relaxations for Problems 1-7}
\label{tab: relaxtion results}
\scalebox{0.9}{
\begin{tabular}{clcccccc}
\hline 
Instance & Relaxation & $LB$ & rgap(\%) & vars. & const. & $K_{exp}$ & Time (s)   \\
\hline
P1 & ECPR \eqref{ECP Relax SGP} & 7.4998 & 87.15 & 14 & 5 & 8 & 0.1  \\
($z^* = 58.38488$) & s-ECPR \eqref{s-ECP Relax SGP} & 56.7598 & \textbf{2.78} & 14 & 11 & 8 & 0.1 \\
&\cite{qu2008global} & 50.7862 & 12.94 & 2 & 2 & - & $< 0.1$ \\ 
&\cite{shen2004global} & 10.2094 & 82.51 & 2 & 2 & - & $< 0.1$ \\ 
&\cite{maranas1997global} & -31.2595 & 153.54 & 4 & 2 & 2 & $< 0.1$ \\
\hline 
P2 & ECPR \eqref{ECP Relax SGP} & 296032.51004 & 36.81 & 19 & 5 & 13 & 0.2 \\
($z^* = 468479.9969$) & s-ECPR \eqref{s-ECP Relax SGP} & 464029.43693 & \textbf{0.95} &  19 &  11 &13 & 0.1 \\
&\cite{qu2008global} & 445430.78105 & 4.92 & 4 & 3 & - & $< 0.1$ \\ 
&\cite{shen2004global} & 281644.36436 & 39.88 & 4 & 3 & - & $< 0.1$ \\ 
&\cite{maranas1997global} &  111263.99926 & 76.25& 7 & 3 & 3 & $< 0.1$\\
\hline
P3 & ECPR \eqref{ECP Relax SGP} & 2.01193 & 49.08 &  40 & 8 & 27 & 0.9 \\
($z^* = 3.95116$) & s-ECPR \eqref{s-ECP Relax SGP} & 3.70697
  & \textbf{6.18} &  40 &  21 &27& 0.8 \\
&\cite{qu2008global} & 3.03844 & 23.10 & 8 &  4 & - &$< 0.1$ \\ 
&\cite{shen2004global} & 2.78004 & 29.64 & 8 &  4 & - & $< 0.1$\\ 
& \cite{maranas1997global} & -4.41660 & 211.78 & 10 & 4 & 2 & 0.1\\ \hline
P4 & ECPR \eqref{ECP Relax SGP} & 2153.54527 & 69.45 & 28 & 10 & 17 & 1.1 \\
($z^* =  7049.24803$) & s-ECPR \eqref{s-ECP Relax SGP} & 6760.93408 & \textbf{4.09} & 28 &  21 & 17 & 1.1\\
&\cite{qu2008global} & 5794.48188 & 17.80 & 8 & 5 & - & $< 0.1$ \\ 
&\cite{shen2004global} & 4718.76663 & 33.06 & 8 & 5 & - & $< 0.1$\\ 
&\cite{maranas1997global} & 1463.42389 & 79.24 & 8 & 5 & - & $< 0.1$\\
\hline
P5 & ECPR \eqref{ECP Relax SGP} &  3471.83273 & 44.16  & 18 & 3 & 13 & 0.1 \\
($z^* =  6217.46549$) & s-ECPR \eqref{s-ECP Relax SGP} &  6019.75009 & \textbf{3.18}  & 18 &  8 & 13 & 0.1 \\
&\cite{qu2008global} & 5537.89651 & 10.93 & 3 & 1 & - & $< 0.1$ \\ 
&\cite{shen2004global} & 5833.22612 & 21.12 & 3 & 1 & - & $< 0.1$\\ 
&\cite{maranas1997global} & 1922.44032 & 69.08 & 8 & 5 & 3 & $< 0.1$\\
\hline
P6 & ECPR \eqref{ECP Relax SGP} & 4139.23598 & 59.11 & 43 & 16 & 28 & 1.9 \\
($z^* =  10122.85643$) & s-ECPR \eqref{s-ECP Relax SGP} & 9865.73588 & \textbf{2.54}  & 43 &  31 & 28 & 1.6\\
&\cite{qu2008global} & 7591.13004 & 25.01 & 5 & 6 & - & $< 0.1$ \\ 
&\cite{shen2004global} & 6614.27439 & 34.66 & 5 & 6 & - & $< 0.1$\\
&\cite{maranas1997global} & 5038.14564 & 50.23 & 11 & 6 & 1 & $< 0.1$\\
\hline
P7 & ECPR \eqref{ECP Relax SGP} & -31.364723 & 62.51 & 21 & 5 & 13 & 0.3 \\
($z^* =  -83.66157$) & s-ECPR \eqref{s-ECP Relax SGP} & --75.54639 & \textbf{9.70}  & 21 &  13 & 13 & 0.2 \\
&\cite{qu2008global} &-67.22207  & 19.65 & 3 & 1 & - & $< 0.1$ \\ 
&\cite{shen2004global} & -58.01929 & 30.65 & 3 & 1 & - & $< 0.1$\\
&\cite{maranas1997global} & -12.76675 & 84.74 & 4 & 3 & 2 & $< 0.1$\\
\hline
\multicolumn{8}{l}{rgap: root gap \qquad  vars.: number of variables \qquad const.: number of constraints}   \\
\multicolumn{8}{l}{$K_{\text{exp}}$: number of exponential conic constraints} \\ \hline
\end{tabular}}
\end{table}
The results for all test problems can be found in Table \ref{tab: relaxtion results}, where for each instance, $z^*$ refers to the optimal solution of the true problem, column $LB$ corresponds to the optimal solution of relaxations (a lower bound for the corresponding $z^*$) and columns $\text{rgap}$ correspond to the relaxation gap percentage obtained from $\frac{100(z^*-LB)}{LB}$. The following observations are obtained from Table \ref{tab: relaxtion results}:
\begin{itemize}
\item s-ECPR provides significantly smaller gaps, demonstrating that it is substantially tighter for the SGP problems considered. Our proposed valid inequalities provide significant advantage, as demonstrated.
\item Linear relaxations \cite{shen2004global} and \cite{qu2008global} have the same number of variables and constraints as the original SGP problems. However, our proposed relaxations involve additional variables and constraints as well as (nonlinear) cones.  This results in long per-relaxation solve times for larger instances such as P4 or P6. Note, however, that the size of the proposed relaxation increases only linearly with respect to the original formulation.
\end{itemize}
\subsection{Sequential Convergence}
In the previous section we demonstrated that our proposed relaxation is substantially stronger, but takes more time to solve compared to linear relaxations. Such a tradeoff may be favorable when in SGP algorithms where relaxation strength governs the total number of iterations involved. This section compares the performance of the proposed sequential ECP optimization algorithm with that of \cite{xu2014global}. For illustration, consider the following problem:

\noindent \textbf{Problem 8}:
\begin{align*}
(\text{P8}) \quad \min \quad & x_1 + x_2 + x_3 \\
& 1 \le x_1x_2 + x_1x_3, \\
& 0.5 \le x_1, x_2, x_3 \le 10.
\end{align*}
Using the bounds on $x$, 
we derive the following s-ECP relaxation:

\begin{align} \label{ex 0 relax}
\min_x \quad & x_1 + x_2 + x_3 \nonumber\\
& 1 \le \gamma_{11} + \gamma_{12} \nonumber, \\
& \tilde{\gamma}_{11} \le \tilde{x}_{1} + \tilde{x}_2 \nonumber, \\
& \tilde{\gamma}_{12} \le \tilde{x}_1 + \tilde{x}_3 \nonumber,\\
& \text{e}^{\tilde{x}_i} \le x_i \qquad  i =1,2,3 \nonumber,\\
& \text{e}^{\tilde{\gamma}_{1j}} \le \gamma_{1j} \qquad j =1,2 \nonumber,\\
& {x}_{i} \le \frac{9.5}{\log(20)}\Big(\tilde{x}_{i} - \log(0.5)\Big) + 0.5 \qquad \quad i =1,2,3 \nonumber,\\
& {\gamma}_{1j} \le \frac{99.75}{\log(400)}\Big(\tilde{\gamma}_{1j}- \log(0.25)\Big) + 0.25, \nonumber \\
& 0.5 \le x_i \le 10, \quad \log(0.1) \le \tilde{x}_i \le \log(10) \qquad i =1,2,3 \nonumber, \\
& 0.25 \le \gamma_{1j} \le 100, \quad \log(0.25) \le \tilde{\gamma}_{1j} \le \log(100) \qquad j =1,2 \nonumber.
\end{align}

Solving this yields an initial solution ${x^{(0)}} = (0.5, 0.5, 0.5)$, $\tilde{x}^{(0)} = (-0.69, -0.69, -0.69)$, and $\tilde{\gamma}^{(0)} = (-1.39, -1.39)$. Applying this to ECP subproblem \eqref{ECP sub SGN}, we achieve the optimal solution $x^* = (1, 0.5, 0.5)$ after 4 iterations.

Now we shall compare this to the sequential GP method of \cite{xu2014global}. We test various initializations. Table \ref{tab:prem results} provides a summary of the results. The first row corresponds to the outcomes of sequential ECP with the initial solution $x^{(0)}$ obtained from the convex hull relaxation of the presented SGP. Subsequent rows display the results of \cite{xu2014global} with different starting points $x^{(0)}$, as specified below.

\begin{table}[htbp]
\caption{Convergence Results for Problem 8}
\scalebox{1}{
\begin{tabular}{c|cccccccc}
\hline
        \textbf{Method} & \textbf{$x^0$} & \textbf{Feasible?} & \textbf{$z^0$} & \textbf{$x^*$} & \textbf{$z^*$} & \textbf{Gap} & \textbf{Iter.} & \textbf{Time (s)} \\
        \hline
        Our work & (0.5, 0.5, 0.5) & No & 1.5 & (1,0.5,0.5) & 2 & 0.25  & 4 & 0.7 \\
        \hline
        \multirow{6}{*}{\cite{xu2014global}} 
        & (1.0, 1.0, 1.0) & Yes & 3.0 & (1, 0.5, 0.5) & 2 & 0.50 & 5 & 0.5 \\
        & (2.0, 2.0, 0.5) & Yes & 4.5 & (1, 0.5, 0.5) & 2 & 1.25 & 5 & 0.7 \\
        & (0.5, 1.0, 0.5) & No & 2 & (1, 0.5, 0.5)& 2 & 0.0 & 6 & 0.5 \\
        & (0.5, 1.4, 0.5) & No & 3.4 & (1, 0.5, 0.5) & 2 & 0.70 & 7 & 0.9 \\
        & (0.5, 5, 10) & Yes & 15.5 & (1, 0.5, 0.5) & 2 & 6.75 & 7 & 0.8 \\
        & (10, 10, 10) & Yes & 30 & (1, 0.5, 0.5) & 2 & 14 & 8 & 0.8 \\
        \hline
    \end{tabular}}
    \label{tab:prem results}
\end{table}
In Table \ref{tab:prem results}, the "Gap" column is defined as $|\frac{z^0 - z^*}{z^*}|$, capturing the difference of the initial solution and the optimal solution in each instance. In the case of sequential ECP, the algorithm converges after a few iterations. 
In contrast, in \cite{xu2014global}, the convergence of the algorithm depends on the quality of the predetermined starting solution $x^{(0)}$, and infeasible starting points can adversely affect the algorithm's performance, requiring additional iterations for convergence \citep{xu2014global}. 

Obtaining a feasible solution for the sequential GP is not always straightforward, suggesting that a readily available feasible solution may not always be at hand for the sequential GP to commence. Thus, sequential GP may require more iterations if initialized with infeasible solutions. To explore this phenomenon, we conducted numerical experiments by generating 100 initial solutions $x^{(0)}$ within the range $0.5 \le x \le 100$ for problem 8 and recorded the performance statistics of the sequential GP in Table \ref{tab:my_label}.

\begin{table}[h]
    \centering
    \caption{Statistics for sequential GP of \cite{xu2014global} for 100 instances of initial solutions}
    \begin{tabular}{c|ccc|ccc}
    \hline
 \% of feasible   & \multicolumn{3}{c|}{Iterations} & \multicolumn{3}{c}{Time (s)} \\
 solutions & max & min & average & max& min  & average \\ \hline
         12 & 9 & 3 & 6.7 & 0.1 & 1.3 & 0.88\\ \hline
    \end{tabular}
    
    \label{tab:my_label}
\end{table}
It can be observed that  $88\%$ of experiments start with infeasible initial solutions, resulting in 6.7 iterations on average to solve SGP problem 8. 

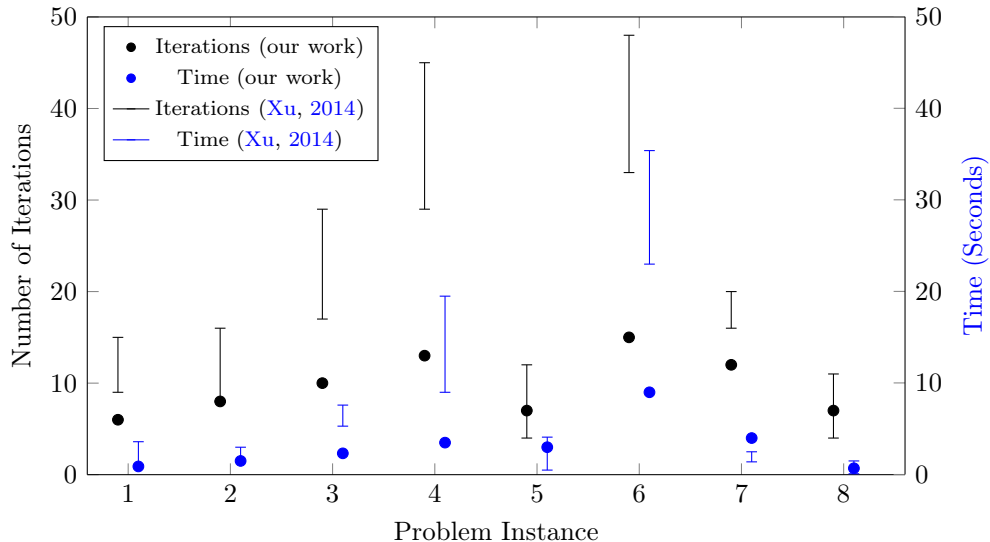
\begin{figure}[t]
\centering
\medskip
\begin{tikzpicture}[font=\small]
\pgfplotsset{width=0.75\textwidth, height=0.35\textheight}
\begin{axis}[
    legend image post style={scale=0.8},
    legend style={{font=\fontsize{8}{8}\selectfont}, at={(0.02,0.98)}, anchor=north west},
    xmin = 0.6,
    xmax = 8.6,
    xtick = {1,2,...,8},
    xlabel = Problem Instance,
    ymin = 0,
    ymax = 50,
    ytick = {0,10,...,50},
    ylabel = {Number of Iterations}, 
    axis y line*=left,
]
\addplot[scatter/classes={
a={mark=*,black}},
scatter,only marks,
scatter src=explicit symbolic] coordinates {(0.9,6)[a] (1.9,8)[a] (2.9,10)[a] (3.9,13)[a]  (4.9,7)[a] (5.9,15)[a] (6.9,12)[a] (7.9,7)[a]};
\addplot[scatter/classes={
a={mark=*,blue}},
scatter,only marks,
scatter src=explicit symbolic] coordinates {(1.1,0.9)[a] (2.1,1.5)[a] (3.1,2.33)[a] (4.1,3.5)[a]  (5.1,3)[a] (6.1,9)[a] (7.1,4)[a] (8.1,0.7)[a]};

\addplot [mark=-, black] plot coordinates {  (0.9,9) (0.9,15) };
\addplot [mark=-, blue] plot coordinates {  (1.1,1.1) (1.1,3.6) };
\addplot [mark=-, black] plot coordinates {  (1.9,8) (1.9,16) };
\addplot [mark=-, blue] plot coordinates {  (2.1,1.4) (2.1,3) };
\addplot [mark=-, black] plot coordinates {  (2.9,17) (2.9,29) };
\addplot [mark=-, blue] plot coordinates {  (3.1,5.3) (3.1,7.6) };
\addplot [mark=-, black] plot coordinates {  (3.9,29) (3.9,45) };
\addplot [mark=-, blue] plot coordinates {   (4.1,9) (4.1,19.5) };
\addplot [mark=-, black] plot coordinates {  (4.9,4) (4.9,12) };
\addplot [mark=-, blue] plot coordinates {  (5.1,0.5) (5.1,4.1) };
\addplot [mark=-, black] plot coordinates {  (5.9,33) (5.9,48) };
\addplot [mark=-, blue] plot coordinates {  (6.1,23) (6.1,35.4) };
\addplot [mark=-, black] plot coordinates {  (6.9,16) (6.9,20) };
\addplot [mark=-, blue] plot coordinates {  (7.1,1.4) (7.1,2.5) };
\addplot [mark=-, black] plot coordinates {  (7.9,4) (7.9,11) };
\addplot [mark=-, blue] plot coordinates {  (8.1,0.1) (8.1,1.5) };

\legend{Iterations (our work), Time (our work), Iterations \citep{xu2014global}, Time \citep{xu2014global}}
\end{axis}

\begin{axis}[
    axis y line*=right,
    ymin=0,
    ymax=50, 
    ytick={0,10,...,50}, 
    ylabel={\color{blue} Time (Seconds)},
    axis x line=none, 
]
\addplot[draw=none] coordinates {(0,0) (1,24)};
\end{axis}
\end{tikzpicture}
\caption{Comparison of the iteration count, on the left axis, and solution times, on the right axis, between \cite{xu2014global}'s approach and our proposed approach.}
\label{simulationfinal}
\end{figure}

To further expand our numerical experiment for other instances, we generated 10 random initial solutions $x^{(0)}$ within the variable bounds to initialize the sequential algorithm proposed by \cite{xu2014global}. Our proposed algorithm does not require a user-defined initialization, and uses the solution of the convex hull relaxation. Figure \ref{simulationfinal} presents the corresponding the iteration count, on the left axis, and solution times, on the right axis, for both approaches. The variability over all 10 instances for \cite{xu2014global}'s approach, depending on the initial solutions $x^{(0)}$, is indicated by the bars. Conversely, our approach relies solely on the relaxation solution, representing a single point for each instance. Figure \ref{simulationfinal} reveals that the proposed algorithm exhibits higher efficiency in the conducted problem instances, showcasing a reduction in both the number of iterations and solution time. This superiority can be attributed primarily to the quality of the proposed ECP relaxations, which offer tight initial convex relaxations to bootstrap the algorithm.


\section{Conclusion} \label{sec: conclusion}
This paper has introduced a novel convex relaxation for Signomial Geometric Programming (SGP) based on Exponential Conic Programming (ECP), setting our work apart from existing relaxations by not requiring bounded variables in the base model. Additionally, we developed a flexible method of strengthening the relaxation, exploiting any bounded variables or monomial terms available to enhance the strength of relaxations by means of valid linear inequalities.

Through computational experiments, we substantiated the effectiveness of our relaxation approach, demonstrating relatively small optimality gaps. Furthermore, we integrated the proposed relaxation into a sequential optimization algorithm designed to solve SGP instances. Compared with a GP method, we demonstrated fewer iterations and overall time needed to convergence. Our method is further distinguished in that it does not rely on a user-provided starting point.

For future work, it is worth exploring the implementation of our proposed relaxations within a branch and cut algorithm, incorporating recent developments in dynamic inequality generation for SGP \citep{xu2022cutting}. This could leverage bound-tightening strategies in preprocessing to enhance the solution procedure, potentially further improving the efficiency and applicability of our methodology. 

\section*{Declaration of Generative AI and AI-assisted technologies in the writing process}
Statement: During the preparation of this work the author(s) used no AI-assisted tools.


\bibliographystyle{elsarticle-harv} 
\bibliography{references}

\end{document}